\documentclass[12pt]{amsart}
\usepackage{amssymb, hyperref}
\usepackage[margin=1in]{geometry}

\newcommand{\Bgp}{{\Z^\N}}

\long\def\forget#1\forgotten{}
\newcommand{\issuenumber}{37}
\newcommand{\issuemonth}{Dec}
\newcommand{\issueyear}{2014}

\newcommand{\R}{\mathbb{R}}

\newcommand{\Op}{\mathrm{O}}

\newcommand{\alephes}{{\aleph_0}}
\newcommand{\ed}{

\forget

\section{Unsolved problems from earlier issues}

\begin{issue}Is $\binom{\Omega}{\Gamma}=\binom{\Omega}{\Tau}$?\end{issue}
\begin{issue}Is $\ufin(\Op,\Omega)=\sfin(\Gamma,\Omega)$?And if not, does $\ufin(\Op,\Gamma)$ imply
$\sfin(\Gamma,\Omega)$?\end{issue}
\stepcounter{issue}
\begin{issue}Does $\sone(\Omega,\Tau)$ imply $\ufin(\Gamma,\Gamma)$?\end{issue}
\stepcounter{issue}
\begin{issue}Is there, in ZFC, an uncountable set satisfying $\sfin(\cB,\cB)$?\end{issue}
\stepcounter{issue}
\begin{issue}Does $X \nin \NON(\cM)$ and $Y\nin\mathsf{D}$ imply that $X\cup Y\nin \COF(\cM)$?\end{issue}
\begin{issue}[CH]Is $\split(\Lambda,\Lambda)$ preserved under finite unions?\end{issue}
\begin{issue}Is $\cov(\cM)=\fo$? (See the definition of $\fo$ in that issue.)\end{issue}
\stepcounter{issue}
\begin{issue}Could there be a Baire metric space $M$ of weight $\aleph_1$ and a partition
$\mathcal{U}$ of $M$ into $\aleph_1$ meager sets where for each ${\mathcal U}'\subset\mathcal U$,
$\bigcup {\mathcal U}'$ has the Baire property in $M$?\end{issue}
\stepcounter{issue} 
\begin{issue}Is there, in ZFC, a set of reals $X$ of cardinality $\fd$ such that all
finite powers of $X$ have Menger's property $\sfin(\Op,\Op)$?\end{issue}
\stepcounter{issue}
\begin{issue}[MA]Is there a set $X\sbst\bbR$, of cardinality continuum, satisfying $\sone(\BO,\BG)$?\end{issue}
\begin{issue}[CH]Is there a totally imperfect $X$ satisfying $\ufin(\Op,\Gamma)$
that can be mapped continuously onto $\Cantor$?\end{issue}
\begin{issue}[CH]Is there a Hurewicz $X$ such that $X^2$ is Menger but not Hurewicz?\end{issue}
\begin{issue}Does the Pytkeev property of $C_p(X)$ imply that $X$ has Menger's property?\end{issue}
\begin{issue}Does every hereditarily Hurewicz space satisfy $\sone(\BG,\BG)$?\end{issue}
\begin{issue}[CH]Is there a Rothberger-bounded $G\le\Bgp$ such that $G^2$ is not Menger-bounded?\end{issue}
\begin{issue}Let $\cW$ be the van der Waerden ideal. Are $\cW$-ultrafilters closed under products?\end{issue}
\begin{issue}Is the $\delta$-property equivalent to the $\gamma$-property $\binom{\Omega}{\Gamma}$?\end{issue}
\stepcounter{issue}\stepcounter{issue}

\forgotten

\general\end{document}}

\newcommand{\Cantor}{{\{0,1\}^\N}}

\newcommand{\fd}{\mathfrak{d}}

\newcommand{\NON}{{\mathsf   {NON}}}\newcommand{\COF}{{\mathsf   {COF}}}

\newcommand{\cM}{\mathcal{M}}
\newcommand{\cov}{\mathsf{cov}}
\newcommand{\cof}{\mathsf{cof}}
\newcommand{\bbR}{\mathbb{R}}
\newcommand{\fo}{\mathfrak{od}}

\renewcommand{\split}{\mathsf{Split}}\newcommand{\bq}{\begin{quote}}\newcommand{\eq}{\end{quote}}
\newcommand{\cB}{\mathcal{B}}\newcommand{\BG}{\cB_\Gamma}
\newcommand{\BO}{\cB_\Omega}

\newcommand{\sone}{\mathsf{S}_1}\newcommand{\sfin}{\mathsf{S}_\mathrm{fin}}
\newcommand{\ufin}{\mathsf{U}_\mathrm{fin}} 
\newcommand{\nin}{\not\in}
\newcommand{\cW}{\mathcal{W}}

\newcommand{\N}{\mathbb{N}}\newcommand{\Z}{\mathbb{Z}}
\newcommand{\sm}{\setminus}\newcommand{\sbst}{\subseteq}
\newcommand{\by}[2]{\par\hfill\emph{#1}, #2}\newcommand{\Tau}{\mathrm{T}}
\newcommand{\CE}{\textsc{CE}}
\newtheorem{thm}{Theorem}[section]\newcommand{\bthm}{\begin{thm}} \newcommand{\ethm}{\end{thm}}
\newtheorem{prop}[thm]{Proposition}\newcommand{\bprp}{\begin{prop}} \newcommand{\eprp}{\end{prop}}
\newtheorem{fact}[thm]{Fact}\newcommand{\bfct}{\begin{fact}} \newcommand{\efct}{\end{fact}}
\newtheorem{prob}[thm]{Problem}\newcommand{\bprb}{\begin{prob}} \newcommand{\eprb}{\end{prob}}
\newtheorem{lem}[thm]{Lemma}\newcommand{\blem}{\begin{lem}} \newcommand{\elem}{\end{lem}}
\newtheorem{claim}[thm]{Claim}\newcommand{\bclm}{\begin{claim}} \newcommand{\eclm}{\end{claim}}
\newtheorem{cor}[thm]{Corollary}\newcommand{\bcor}{\begin{cor}} \newcommand{\ecor}{\end{cor}}
\newtheorem{conj}[thm]{Conjecture}\newcommand{\bcnj}{\begin{conj}} \newcommand{\ecnj}{\end{conj}}
\theoremstyle{definition}\newtheorem{defn}[thm]{Definition}\newcommand{\bdfn}{\begin{defn}} \newcommand{\edfn}{\end{defn}}
\theoremstyle{remark}\newtheorem{rem}[thm]{Remark}\newcommand{\brem}{\begin{rem}} \newcommand{\erem}{\end{rem}}
\newtheorem{cnv}[thm]{Convention}\newcommand{\bcnv}{\begin{cnv}} \newcommand{\ecnv}{\end{cnv}}
\newtheorem{exam}[thm]{Example}\newcommand{\bexm}{\begin{exam}} \newcommand{\eexm}{\end{exam}}
\newtheorem{issue}{Issue}\newcommand{\bpf}{\begin{proof}} \newcommand{\epf}{\end{proof}}
\newcommand{\be}{\begin{enumerate}}\newcommand{\ee}{\end{enumerate}}\newcommand{\bi}{\begin{itemize}}
\newcommand{\ei}{\end{itemize}}
\newcommand{\general}{\small\vfill\par\noindent\hrulefill\par
\noindent\textbf{Previous issues.} 
\url{http://front.math.ucdavis.edu/search?\&t=\%22SPM+Bulletin\%22}
\\[0.1cm]
\textbf{Contributions and free subscription.} Email \url{tsaban@math.biu.ac.il}.
}

\newcommand{\link}[1]{\par\hfill{\url{#1}}}

\newcommand{\arXivl}[4]{\subsection{#2}{#4}\par\hfill{\arx{#1}}\par\hfill\emph{#3}}
\newcommand{\arXiv}[3]{\subsection{#2}\mbox{}\par\hfill{\arx{#1}}\par\hfill\emph{#3}}

\newcommand{\AMS}[3]{\subsection{#1}\mbox{}\par\hfill{\url{#3}}\par\hfill\emph{#2}}

\newcommand{\arx}[1]{\url{http://arxiv.org/abs/#1}}

\title[$\mathcal{SPM}$ Bulletin \textbf{\issuenumber} (\issuemonth{} \issueyear)]{%
$\mathcal{SPM}$ Bulletin\\[0.5cm]
Issue number \issuenumber: \issuemonth{} \issueyear{} \CE{}}

\begin{document}
\maketitle


\section{Editor's note}

In the year 2014, the field of selection principles found its way into several additional, 
fascinating mathematical realms. The field enters the concensus
as a mainstream part of set theory and topology, and as a promising direction for
young researchers that are well trained in these fields.

Some of the interesting developments in the field during 2014 are reported in this issue.

\medskip

With best regards,

\by{Boaz Tsaban}{tsaban@math.biu.ac.il}

\hfill \texttt{http://www.cs.biu.ac.il/\~{}tsaban}

\section{Long announcements}

\arXivl{1312.6081}
{A note on condensations of function spaces onto $\sigma$-compact and
  analytic spaces}
{Miko{\l}aj Krupski}
{Modifying a construction of W. Marciszewski we prove (in ZFC) that there
exists a subspace of the real line $\mathbb{R}$, such that the realcompact
space $C_p(X)$ of continuous real-valued functions on $X$ with the pointwise
convergence topology does not admit a continuous bijection onto a
$\sigma$-compact space. This answers a question of Arhangel'skii.}

\arXivl{1401.2283}
{Mathias forcing and combinatorial covering properties of filters}
{David Chodounsk\'y, Du\v{s}an Repov\v{s}, Lyubomyr Zdomskyy}
{We give topological characterizations of filters $F$ on $w$ such that the
Mathias forcing $M_F$ adds no dominating reals or preserves ground model
unbounded families. This allows us to answer some questions of Brendle,
Guzm\'an, Hru\v{s}\'ak, Mart\'{\i}nez, Minami, and Tsaban.}

\arXivl{1401.6061}
{Baire spaces and infinite games}
{Fred Galvin and Marion Scheepers}
{It is well known that if the nonempty player of the Banach--Mazur game has a
winning strategy on a space, then that space is Baire in all powers even in the
box topology. The converse of this implication may be true also: We know of no
consistency result to the contrary, and in this paper establish the consistency
of the converse relative to the consistency of the existence of a proper class
of measurable cardinals.}

\arXivl{1405.4929}
{Selective Games on Binary Relations}
{Rodrigo R. Dias and Marion Scheepers}
{We present a unified approach, based on dominating families in binary
relations, for the study of topological properties defined in terms of
selection principles and the games associated to them.}

\arXivl{1405.5568}
{Some observations on filters with properties defined by open covers}
{Rodrigo Hern\'andez--Guti\'errez, Paul J. Szeptycki}
{We study the relation between the Hurewicz and Menger properties of filters
considered topologically as subspaces of $P(\omega)$ with the Cantor set
topology.}

\arXivl{1405.7208}
{Combinatorial aspects of selective star covering properties in
  $\Psi$-spaces}
{Boaz Tsaban}
{Which Isbell--Mr\'owka spaces ($\Psi$-spaces) satisfy the star version of
Menger's and Hurewicz's covering properties? Following Bonanzinga and Matveev,
this question is considered here from a combinatorial point of view. An example
of a $\Psi$-space that is (strongly) star-Menger but not star-Hurewicz is
obtained. The PCF-theory function $\kappa\mapsto\cof([\kappa]^\alephes)$ is a
key tool. Using the method of forcing, a complete answer to a question of
Bonanzinga and Matveev is provided.

  The results also apply to the mentioned covering properties in the realm of
Pixley--Roy spaces, to the extent of spaces with these properties, and to the
character of free abelian topological groups over hemicompact $k$ spaces.}

\arXivl{1406.0692}
{Countable dense homogeneous filters and the Menger covering property}
{Du\v{s}an Repov\v{s}, Lyubomyr Zdomskyy, and Shuguo Zhang}
{In this note we present a ZFC construction of a non-meager filter which fails
to be countable dense homogeneous. This answers a question of
Hern\'andez-Guti\'errez and Hru\v{s}\'ak. The method of the proof also allows
us to obtain a metrizable Baire topological group which is strongly locally
homogeneous but not countable dense homogeneous.

Journal reference: Fundamenta Mathematicae 224 (2014), 233--240.
}

\arXivl{1406.0696}
{Productively Lindel\"of spaces and the covering property of Hurewicz}
{Du\v{s}an Repov\v{s} and Lyubomyr Zdomskyy}
{We prove that under certain set-theoretic assumptions every productively
Lindel\"of space has the Hurewicz covering property, thus improving upon some
earlier results of Aurichi and Tall.

Journal reference: Topology and its Applications 169 (2014), 16--20.
}

\arXivl{1409.7337}
{Point Networks for Special Subspaces of $\R^{\kappa}$}
{Ziqin Feng and Paul Gartside}
{Uniform characterizations of certain special subspaces of products of lines
are presented. The characterizations all involve a collection of subsets (base,
almost subbase, network or point network) organized by a directed set. New
characterizations of Eberlein, Talagrand and Gulko compacta follow.}

\arXivl{1406.3062}
{Luzin and Sierpi\'nski sets, some nonmeasurable subsets of the plane and
  additive properties on the line}
{Marcin Michalski, Szymon \.Zeberski}
{In this paper we shall introduce some nonmeasurable and completely
nonmeasurable subsets of the plane with various additional properties, e.g.,
being Hamel basis, intersection with each line is a super Luzin / Sierpi\'nski
set. Also some additive properties of Luzin and Sierpi\'nski sets on the line
are investigated. One of the result is that for a Luzin set $L$ and a
Sierpi\'nski set $S$, the set $L+S$ cannot be a Bernstein set.}

\arXivl{1406.7805}
{Between countably compact and $\omega$-bounded}
{Istv\'an Juh\'asz and Lajos Soukup and Zolt\'an Szentmikl\'ossy}
{Given a property $P$ of subspaces of a $T_1$ space $X$, we say that $X$ is
{\em $P$-bounded} iff every subspace of $X$ with property $P$ has compact
closure in $X$. Here we study $P$-bounded spaces for the properties $P \in
\{\omega D, \omega N, C_2 \}$ where $\omega D \, \equiv$ "countable discrete",
$\omega N \, \equiv$ "countable nowhere dense", and $C_2 \,\equiv$ "second
countable". Clearly, for each of these $P$-bounded is between countably compact
and $\omega$-bounded.
  We give examples in ZFC that separate all these boundedness properties and
their appropriate combinations. Consistent separating examples with better
properties (such as: smaller cardinality or weight, local compactness, first
countability) are also produced.
  We have interesting results concerning $\omega D$-bounded spaces which show
that $\omega D$-boundedness is much stronger than countable compactness:
\be
\item Regular $\omega D$-bounded spaces of Lindel\"of degree $<
cov(\mathcal{M})$ are $\omega$-bounded.
\item Regular $\omega D$-bounded spaces of countable tightness are
$\omega N$-bounded, and if $\mathfrak{b} > \omega_1$ then even
$\omega$-bounded.
\item If a product of Hausdorff space is $\omega D$-bounded then all but
one of its factors must be $\omega$-bounded.
\item Any product of at most $\mathfrak{t}$ many Hausdorff $\omega
D$-bounded spaces is countably compact.
\ee
  As a byproduct we obtain that regular, countably tight, and countably compact
spaces are discretely generated.}

\arXivl{1407.7495}
{When is a space Menger at infinity?}
{Leandro F. Aurichi, Angelo Bella}
{We try to characterize those Tychonoff spaces $X$ such that $\beta X
\setminus X$ has the Menger property.}

\arXivl{407.7437}
{Algebra in the Stone--\v{C}ech compactification,
selections, and additive combinatorics}
{Boaz Tsaban}
{The algebraic structure of the  Stone--\v{C}ech compactification of a semigroup,
and methods from the theory of selection principles, are used to establish qualitative coloring theorems
extending the Milliken--Taylor Theorem and, consequently, Hindman's Finite Sums Theorem.
The main result is the following one (definitions provided in the main text):
Let $X$ be a Menger space, and $\mathcal{U}$ be a point-infinite open cover of $X$ with no finite subcover.
Consider the complete graph, whose vertices are the open sets in $X$.
For each finite coloring of the edges of this graph,
there are disjoint finite subsets $\mathcal{F}_1,\mathcal{F}_2,\dots$ of the cover $\mathcal{U}$
whose unions $V_1 := \bigcup\mathcal{F}_1, V_2 := \bigcup\mathcal{F}_2,\dots$
have the following properties:
\begin{enumerate}
\item The family $\{V_1,V_2,\dots\}$ is a point-infinite cover of $X$.
\item The sets $\bigcup_{n\in F}V_n$ and $\bigcup_{n\in H}V_n$ are distinct for all nonempty finite sets
$F<H$.
\item All edges $\bigl\{\,\bigcup_{n\in F}V_n, \bigcup_{n\in H}V_n\,\bigr\}$, for nonempty
finite sets $F<H$, have the same color.
\end{enumerate}
A purely combinatorial consequence of this result is provided.
A self-contained introduction to the necessary parts of the needed theories,
modulo the definition and elementary properties of ultrafilters, is provided.}

\arXivl{1409.0579}
{CH and the Moore--Mrowka Problem}
{Alan Dow and Todd Eisworth}
{We show that the Continuum Hypothesis is consistent with all regular spaces
of hereditarily countable $\pi$-character being C-closed. This gives us a model
of ZFC in which the Continuum Hypothesis holds and compact Hausdorff spaces of
countable tightness are sequential.}

\arXivl{1412.1494}
{On $\mathfrak{P}$-spaces and related concepts}
{S. S. Gabriyelyan, J. Kakol}
{The concept of the strong Pytkeev property, recently introduced by Tsaban and
Zdomskyy in [32], was successfully applied to the study of the space $C_c(X)$
of all continuous real-valued functions with the compact-open topology on some
classes of topological spaces $X$ including \v{C}ech-complete Lindel\"{o}f
spaces. Being motivated also by several results providing various concepts of
networks we introduce the class of $\mathfrak{P}$-spaces strictly included in
the class of $\aleph$-spaces. This class of generalized metric spaces is closed
under taking subspaces, topological sums and countable products and any space
from this class has countable tightness. Every $\mathfrak{P}$-space $X$ has the
strong Pytkeev property. The main result of the present paper states that if
$X$ is an $\aleph_0$-space and $Y$ is a $\mathfrak{P}$-space, then the function
space $C_c(X,Y)$ has the strong Pytkeev property. This implies that for a
separable metrizable space $X$ and a metrizable topological group $G$ the space
$C_c(X,G)$ is metrizable if and only if it is Fr\'{e}chet-Urysohn. We show that
a locally precompact group $G$ is a $\mathfrak{P}$-space if and only if $G$ is
metrizable.}

\arXivl{1412.4268}
{The strong Pytkeev property in topological spaces}
{Taras Banakh and Arkady Leiderman}
{A topological space $X$ has the strong Pytkeev property at a point $x\in X$
if there exists a countable family $\mathcal N$ of subsets of $X$ such that for
each neighborhood $O_x\subset X$ and subset $A\subset X$ accumulating at $x$,
there is a set $N\in\mathcal N$ such that $N\subset O_x$ and $N\cap A$ is
infinite. We prove that for any $\aleph_0$-space $X$ and any space $Y$ with the
strong Pytkeev property at a point $y\in Y$ the function space $C_k(X,Y)$ has
the strong Pytkeev property at the constant function $X\to \{y\}\subset Y$. If
the space $Y$ is rectifiable, then the function space $C_k(X,Y)$ is rectifiable
and has the strong Pytkeev property at each point. We also prove that for any
pointed spaces $(X_n,*_n)$, $n\in\omega$, with the strong Pytkeev property
their Tychonoff product and their small box-product both have the strong
Pytkeev property at the distinguished point. We prove that a sequential
rectifiable space $X$ has the strong Pytkeev property if and only if $X$ is
metrizable or contains a clopen submetrizable $k_\omega$-subspace. A locally
precompact topological group is metrizable if and only if it contains a dense
subgroup with the strong Pytkeev property.}

\section{Short announcements}\label{RA}

\arXiv{1311.2285}
{Ultrafilter convergence in ordered topological spaces}
{Paolo Lipparini}

\arXiv{1311.2330}
{Cardinal invariants for $\kappa$-box products}
{W. W. Comfort and Ivan S. Gotchev}

\arXiv{1311.6544}
{The non-Urysohn number of a topological space}
{Ivan S. Gotchev}

\arXiv{1312.6946}
{Scattered Subsets of Groups}
{T.O. Banakh, I.V. Protasov, S.V. Slobodianiuk}

\arXiv{1312.7527}
{$\mathcal I^{\mathcal K}$-Cauchy functions}
{Pratulananda Das, Martin Sleziak and Vladim\'ir Toma}

\arXiv{1401.2289}
{Metrizable images of the Sorgenfrey line}
{Mikhail Patrakeev}

\arXiv{1401.3132}
{On subcontinua and continuous images of $\beta\R\sm\R$}
{Alan Dow and Klaas Pieter Hart}

\arXiv{1402.1589}
{Compact spaces, lattices, and absoluteness: a survey}
{Wies{\l}aw Kubi\'s}

\arXiv{1403.1909}
{On Haar meager sets}
{U. B. Darji}

\arXiv{1404.3703}
{Coherent ultrafilters and nonhomogeneity}
{Jan Star\'y}

\arXiv{1405.2151}
{The covering number of the difference sets in partitions of $G$-spaces
  and groups}
{Taras Banakh and Mikolaj Fraczyk}

\arXiv{1405.7899}
{Between Polish and completely Baire}
{Andrea Medini, Lyubomyr Zdomskyy}

\arXiv{1406.1405}
{Partitions of $2^{\omega}$ and completely ultrametrizable spaces}
{William R. Brian and Arnold W. Miller}

\subsection{Metrization conditions for topological vector spaces with Baire 
type properties }
\mbox{}
\par\hfill{DOI: 10.1016/j.topol.2014.05.007}
\par\hfill\emph{S. S. Gabriyelyan, J. K\c{a}kol}

\arXiv{1409.3922}
{Classifying invariant $\sigma$-ideals with analytic base on good Cantor
  measure spaces}
{Taras Banakh, Robert Ralowski, Szymon Zeberski}

\arXiv{1410.0559}
{A homogeneous space whose complement is rigid}
{Andrea Medini, Jan van Mill, Lyubomyr Zdomskyy}

\arXiv{1409.7902}
{P-domination and Borel sets}
{D. Basile, U. B. Darji}

\arXiv{1410.1504}
{A characterization of Tychonoff spaces with applications to
  paratopological groups}
{Taras Banakh and Alex Ravsky}

\AMS{Discontinuity of multiplication and left translations in $\beta G$}
{Yevhen Zelenyuk}
{www.ams.org/journal-getitem?pii=S0002-9939-2014-12267-2}

\arXiv{1410.6081}
{P-sets and minimal right ideals in $\N^*$}
{William R. Brian}

\arXiv{1411.0633}
{Measure of compactness for filters in product spaces:
  Kuratowski--Mr\`owka in CAP revisited}
{Fr\'ed\'eric Mynard and William Trott}

\arXiv{1412.2216}
{On the $C_k$-stable closure of the class of (separable) metrizable
  spaces}
{T. Banakh, S. Gabriyelyan}

\arXiv{1412.2240}
{The $C_p$-stable closure of the class of separable metrizable spaces}
{T. Banakh, S. Gabriyelyan}

\arXiv{1412.3254}
{A survey on structural Ramsey theory and topological dynamics with the
  Kechris--Pestov--Todorcevic correspondence in mind}
{Lionel Nguyen Van Th\'e}

\AMS{A new class of Ramsey-classification theorems and their applications
in the Tukey theory of ultrafilters, Part 2}
{Natasha Dobrinen; Stevo Todorcevic}
{http://www.ams.org/journal-getitem?pii=S0002-9947-2014-06122-9}

\AMS{On the concept of analytic hardness}
{Janusz Pawlikowski}
{http://www.ams.org/journal-getitem?pii=S0002-9939-2014-12422-1}

\ed